\def\N{{\bf N}}

\def\CC{{\rm\kern.24em
   \vrule width.02em
       height1.4ex depth-.05ex
   \kern-.26em C}}
\def\QQ{{\rm\kern.24em
   \vrule width.02em
       height1.4ex depth-.05ex
   \kern-.26em Q}}
\def\PP{{\rm I\kern-.25em P}}         \def\RR{{\rm I\kern-.25em R}}
\def\DD{{\rm I\kern-.25em D}}         \def\EE{{\rm I\kern-.25em E}}
\def\FF{{\rm I\kern-.25em F}}         \def\NN{{\rm I\kern-.23em N}}
\def\RRp{{\rm I\kern-.25em R}_{+}}
\def\IND{{\rm 1\kern-.25em I}}
\def\GA{{\mbox{\gothic{A}}}}

\documentclass[12pt,a4paper]{article}

\usepackage[T2A]{fontenc}
\usepackage[russian,english]{babel}
\usepackage[cp866]{inputenc}
\font\gothic = eufm10
 \def\casure{\stackrel{\mbox{a.s.}}{\longrightarrow}}
 \newtheorem{theorem}{Theorem}
\newtheorem{lemma}{Lemma}

\author{\MakeUppercase{E.E.\, Permyakova}}

\title{Limit theorems for random processes with random time substitution}
\date{Izv. Vuzov Mathematic, 12.2008}
\begin{document}
\maketitle
\label{firstpage}

  Abstract. 
In this paper  the  sufficient conditions for convergence in Skorokhod space $D[0,1]$ of sequence of random processes with random  time substitution are obtained.


{\bf\Large Introduction. }

In this paper the problem of convergence in distribution for sequence of random processes with random time substitution in Skorokhod space $D[0,1]$  is considered.  
We prove the limit theorem for random processes with random time substitution, some its corollaries and applications to the Insurance Mathematic. 

In the monograph \cite{silvestrov1974} D.S. Silvestrov  studies the convergence of sequences of random processes with random time substitution (limit theorems for complex random functions).
 
In this paper we investigate the convergence of a sequences of superpositions of random processes, when the outer and inner random processes are independent. We require that sequences of outer and inner processes considered converge and impose some additional conditions on them. On the inner processes we impose the conditions of increase and non-negativity of trajectories. It is necessary to ensure the correctness of the definition of superposition of random processes. Also we impose two conditions ($A$) and ($B$) on outer and inner processes relatively. The fulfillment of one of them involves the convergence of the sequence of random processes with random time substitution (Theorem 1). We adduce several examples which show that the convergence is not true when we weaken this  conditions.

\section{Notation and preliminary results}

Denote:
$\Delta[0,k]$ $(k \in \N)$ is a class of strictly increasing continuous mapping of the segment  $[0,k]$ on it self such that
$\lambda(0)=0, \ \lambda(k)=k$; $D[0,k]$ (relatively,
$D[0,\infty)$) is a Skorokhod space, i.e. the space of functions defined in the segment $[0,k]$ (interval $[0,\infty)$)
 in $\RR$, right-continuous and with a finite limit on the left. In the space $D[0,k]$ we will consider the Skorokhod's metric $$\rho_k(x,y)=$$
$$=\inf\{\varepsilon>0: \exists \lambda\in \Delta[0,k], \
\sup_{0\le t\le k}|x(t)-y(\lambda(t))|\le \varepsilon,\ \
\sup_{0\le t\le k}|\lambda(t)-t|\le \varepsilon\},$$ $x,y\in
D[0,k]$. By
$$\rho_{\infty}(x,y)=\sum_{k=1}^{\infty}\frac{1}{2^k}\frac{\rho_k(x,y)}{1+\rho_k(x,y)},\,\,\,\,\,x,y\in
D[0,\infty) $$ we will denote the metric in $D[0,\infty)$. The sufficient conditions of weak convergence in this space for example  are in $\cite{Stone}$.

Let $A$ is a measurable subset of $D[0,k]$ ($D[0,\infty)$). We will say that 
$X$ is a random process in $A$ or that almost all trajectories of random process $X$ belong to  $A$ if there exists a  Borel probability  measure  ${\cal L} (X)$ in $D[0,k]$ ($D[0,\infty)$) such that  ${\cal L} (X)(A)=1$. If $X$
is a random process in $D[0,k]$ then measure ${\cal L} (X)$ is	
uniquely defined by random process $X$  and is called the distribution of random process $X$. Also by ${\cal L} (G)$ we will denote the distribution of random element $G$.
If it does not lead to contradictions we will denote by the same symbol the random process and the random element with overlapping distributions.

Let $X, \ X_n, \ n\in\NN $ are the random processes in $D[0,k]$.
We will write $X_n\stackrel{d}{\to}X$ in $D[0,k]$ at $n\to\infty$,
if the sequence of distributions  ${\cal L} (X_n)$ weak converges to distribution ${\cal L} (X)$ in  $D[0,k]$. By symbol 
$\casure$ we will denote the almost sure convergence.

Let $X'_n$ is a sequence of random processes in
$D[0,\infty)$ and  $\Lambda_n$ is a sequence of random processes with non-decreasing  non-negativity almost sure finite trajectories in $D[0,1]$ such that  $\Lambda_n(0)=0$. Let
$X'_n\stackrel{d}{\to}X'$ at $n\to\infty$ in $D[0,\infty)$,
$\Lambda_n\stackrel{d}{\to}\Lambda$ at $n\to\infty$ in $D[0,1]$,
$X'_n$ and $\Lambda_n$ are independent for all $n\in\NN$. Consider the random processes with random time substitution $X_n(t)=X'_n(\Lambda_n(t)), \
t\in[0,1]$,   $X(t)=X'(\Lambda(t)), \ t\in[0,1]$. Since the trajectories of  $\Lambda_n, \ \Lambda$ are non-decreasing, nonnegative and right-conditions, the trajectories of $X_n, \ X $ are in the space $D[0,1]$.

Consider the space  $E \equiv D[0,\infty)\times B[0,1]$, where $B[0,1]$ denotes the set of non-decreasing non-negative function of space $D[0,1]$. Introduce the metric in  $E$ by equation:  $$\rho_E(z_1,z_2)=\rho_{\infty}(x_1,x_2)+\rho_1(y_1,y_2), \ \ (x_1,y_1), \ (x_2,y_2)\in E.$$
The space $E$ has the Tihonov topology.
Consider the operator $F$ of the space $E$ in space $D[0,1]$: $$F(x,y)=x\circ y, $$ where $\circ$ denotes the superposition of functions. We will interest the subsets of $E$, where operator $F$ is continuous.
Denote:
$E_1=D[0,\infty)\times \Pi [0,1]$, where $\Pi [0,1]\subset B[0,1]$ is a space of strictly increasing continuous functions satisfying the condition:
$$ \mbox{ exists } C>0 \mbox{ such that }  \  f(1)=C \ \mbox{for all } f \in \Pi [0,1]. \eqno(*)$$ Let $E_2=C[0,\infty)\times B[0,1]$. It is easy to see that $E_1, E_2 \subset E$. 
In the spaces $B[0,1]$ and $\Pi[0,1]$ we will consider the metric $\rho_1$.

The following lemmas are needed for further.

\begin{lemma} 
Operator $F: E_1 \rightarrow D[0,1]$ is continuous.
\end{lemma}

{\bf Proof.}
Let the functions $g_n, \ g$ belong to Skorokhod space $D[0,\infty)$, the functions $\gamma_n, \ \gamma$ belong to the space $\Pi [0,1]$ and
$$\rho_{\infty}(g_n,g)\to 0 \mbox{ at } n\to \infty, \ \  \rho_1(\gamma_n,\gamma)\to 0 \mbox{ at } n\to \infty.$$ Show, that $\rho_1(g_n\circ \gamma_n, g\circ\gamma)\to 0$ at $n\to\infty$. 
Let the function $\lambda_n\in \Delta[0,1]$ is arbitrary.
Denote $\mu_n(t)\equiv \gamma_n\lambda_n^{-1}\gamma^{-1}(t).$ It is easy to see that $\mu_n\in\Delta[0,\gamma(1)].$ Then
 $$ \sup_{0\le t\le 1}|g_n(\gamma_n(t))-g(\gamma(\lambda_n(t)))|=$$
 $$= \sup_{0\le t\le 1}|g_n(\mu_n(\gamma(\lambda_n(t)))-g(\gamma(\lambda_n(t)))|= \sup_{0\le t\le \gamma(1)}|g_n(\mu_n(t))-g(t)|. $$
Let $\varepsilon>0$.
The convergence $g_n\to g$ at $n\to\infty$ in $D[0,\infty)$ implies that exists $n_1\in\NN$ such that for all $n>n_1$ there exists a function $\nu_n \in \Delta[0,\gamma(1)]$ satisfying the inequalities:
\begin{equation} \label{lemma1}
\sup_{0\le t\le \gamma(1)}|g_n(\nu_n(t))-g(t)|<\varepsilon, \ \  \sup_{0\le t\le \gamma(1)}|\nu_n(t)-t|<\varepsilon.
\end{equation}
Choose $\lambda_n(t)=\gamma^{-1}\nu_n^{-1}\gamma_n(t).$ It is easy to see that $\lambda_n \in \Delta[0,1]$. 
Note that
$$\sup_{0\le t \le 1}|\lambda_n(t)-t|= \sup_{0\le t \le 1}|\gamma^{-1}\nu_n^{-1}\gamma_n(t)-t|=\sup_{0\le t \le \gamma_n(1)}|\gamma^{-1}(\nu_n^{-1}(t))-\gamma_n^{-1}(t)|\le$$ 
$$\le \sup_{0\le t \le \gamma_n(1)}|\gamma^{-1}(\nu_n^{-1}(t))-\gamma^{-1}(t)|+\sup_{0\le t \le \gamma_n(1)}|\gamma^{-1}(t)-\gamma_n^{-1}(t)|= $$
$$=\sup_{0\le t \le \gamma(1)}|\gamma^{-1}(\nu_n^{-1}(t)-t+t)-\gamma^{-1}(t)|+\sup_{0\le t \le \gamma(1)}|\gamma^{-1}(t)-\gamma_n^{-1}(t)|.$$
Let $\delta>0$ is arbitrary. The inequalities $(\ref{lemma1})$ imply that we can choose $n_2\in\NN$ such that for all $n>n_2$
 $$\sup_{0\le t\le \gamma(1)}|\nu_n(t)-t|<\delta.$$
The continuity of function $\gamma^{-1}(t)$ implies that exists $0<\delta<\varepsilon$ such that the inequality $$\sup_{0\le t\le \gamma(1)}|\nu_n(t)-t|<\delta \mbox{ implies }\sup_{0\le t \le \gamma(1)}|\gamma^{-1}(\nu_n^{-1}(t)-t+t)-\gamma^{-1}(t)|<\frac{\varepsilon}{2} $$ is true.

The convergence $\gamma_n\to\gamma $ at $n\to\infty$ in $D[0,1]$ and continuity of $\gamma(t)$ imply  (see $\cite{bil}$, sect. 14) the convergence  $\gamma_n\to\gamma $ at $n\to\infty$ in the 	
uniform norm of space $C[0,1]$. The monotony and continuity of functions $\gamma(t), \gamma_n(t)$ imply the convergence  $\gamma^{-1}_n\to\gamma^{-1} $ at $n\to\infty$ in $C[0,\gamma(1)]$. Therefore exists $n_3\in\NN$ such that  for all $n>n_3$
$$\sup_{0\le t \le \gamma(1)}|\gamma^{-1}_n(t)-\gamma^{-1}(t)|<\frac{\varepsilon}{2}.$$
Thus for all $n>\max\{n_1,n_2,n_3\}$
$$\sup_{0\le t\le 1}|g_n(\gamma_n(t))-g(\gamma(\lambda_n(t)))|<\varepsilon, \ \ \sup_{0\le t \le 1}|\lambda_n(t)-t|<\varepsilon.$$ By reason of the arbitrary choice of $\varepsilon$ this implies the convergence  $\rho_1(g_n\circ \gamma_n, g\circ\gamma)\to 0$ at $n\to\infty$.
Lemma is proved.

\begin{lemma}
Operator $F: E_2 \rightarrow D[0,1]$ is continuous.
\end{lemma}
{\bf Proof. }
Let the functions $g_n, \ g$ belong to the space  $C[0,\infty)$ and the functions $\gamma_n, \ \gamma$ to the space $B[0,1]$; and the following conditions are satisfied:
$$\rho_{\infty}(g_n,g)\to 0 \mbox{ at } n\to \infty, \ \  \rho_1(\gamma_n,\gamma)\to 0 \mbox{ at } n\to \infty.$$
Show that $$\rho_1(g_n\circ\gamma_n,g\circ\gamma)\to 0 \mbox{ at } n \to\infty. $$ Note that
$$\rho_1(g_n\circ\gamma_n,g\circ\gamma)\le \sup_{0\le t\le 1}|g_n(\gamma_n(t))-g(\gamma_n(t))|+$$
$$+\inf\left\{\beta>0: \exists \lambda\in\Delta[0,1], \  \sup_{0\le t\le 1}|g(\gamma_n(t))-g(\gamma(\lambda(t))|<\beta\right., $$
$$\ \left.\sup_{0\le t\le 1}\left|\lambda(t)-t\right|<\beta\right\}.
\eqno{(1)}
$$ Let $\varepsilon>0$. Consider the right part of inequality (1). We are:
$$\sup_{0\le t\le 1}|g_n(\gamma_n(t))-g(\gamma_n(t))|\le\sup_{0\le t\le\gamma_n(1)}|g_n(t)-g(t)|.   $$
The convergence  $\gamma_n\to\gamma$ at
$n\to\infty$ in $D[0,1]$ implies  the finiteness of sequence
$\gamma_n$ in metric $\rho_1$. Therefore there exists a number $N\in\NN$ such that $\gamma_n(1)\le N$ for all
$n\in\NN$. The continuity of  $g$ implies
 (see. \cite{bil}, s. 14) the convergence of sequence
 $g_n$ in uniform norm. Thus there exists 
$n_1\in\N$ such that for all $n>n_1$
$$\sup_{0\le t\le\gamma_n(1)}|g_n(t)-g(t)|
\le\sup_{0\le t\le
N}|g_n(t)-g(t)|<\frac{\varepsilon}{2}.
\eqno{(2)}$$

The continuity of  $g$ implies that exists
$0<\delta<\frac{\varepsilon}{2}$ such that
$$|g(t)-g(s)|<\frac{\varepsilon}{2} $$ for all $0\le t,s\le N$ and the condition $|t-s|<\delta$ is satisfied.
 Thus $\gamma_n\to \gamma \mbox{ at } n\to\infty \mbox{ in } D[0,1] $ there exists $n_2\in\N$ such that for all $n>n_2$ we can choose a function $\lambda_n\in\Delta[0,1]$ with property:
$$  \sup_{0\le t\le 1}|\gamma_n(t)-\gamma(\lambda_n(t))|<\delta, \ \sup_{0\le t\le 1}\left|\lambda_n(t)-t\right|<\delta.$$
 Therefore for all $n>n_2$
$$\inf\left\{\beta>0: \exists \lambda_n\in\Delta[0,1], \  \sup_{0\le t\le 1}|g(\gamma_n(t))-g(\gamma(\lambda_n(t))|<\beta, \right.$$
$$\ \left.\sup_{0\le t\le 1}|\lambda_n(t)-t|<\beta \right\}<\frac{\varepsilon}{2}.\eqno{(3)}$$
Let $n_0=\max\{n_1, n_2\}$. Then  (1), (2) and (3) imply
$$\rho_1(g_n\circ\gamma_n,g\circ\gamma)<\varepsilon $$ for all $n>n_0$. Lemma is proved.

\section{The main results.}
Further we will consider the sequences of random processes $X'_n$ and $\Lambda_n$ meet one of the following conditions:

$$  \Lambda_n\stackrel{d}{\to}\Lambda\mbox{ at } n\to\infty \mbox{ in } \Pi[0,1], \ X'_n\stackrel{d}{\to}X' \mbox{ at } n\to\infty \mbox{ in }D[0,\infty) \eqno(A)$$

$$\Lambda_n\stackrel{d}{\to}\Lambda \mbox{ at } n\to\infty \mbox{ in
} B[0,1],  \  X'_n\stackrel{d}{\to}X'\mbox{ at } n\to\infty\mbox{ in
}C[0,\infty)   \eqno(B) $$

\begin{theorem}
Let  the random processes $X'_n$ and $\Lambda_n$, $X'$ and $\Lambda$ are independent for all $n\in\NN$ and the condition $(A)$ or $(B)$ is satisfied.  Тhen
$$X_n\stackrel{d}{\to}X \mbox{ at } n\to\infty 
\mbox{ in } D[0,1].$$
\end{theorem}
{\bf Proof.}
Let the condition $(A)$ is satisfied.
By Skorokhod theorem about one probability space (see, for example, theorem 11 in section V in \cite{Bul}), it exists the probability space $(\Omega'_1, \GA_1, \PP_1)$ and the random processes $G'_n:\Omega'_1\rightarrow D[0,\infty)$ such that

\medskip

$1) {\cal L} (X'_n)={\cal L} (G'_n); \ \ $

\medskip

$ 2) \ G'_n\casure G' \mbox{ at } n\to\infty \mbox{ in }
D[0,\infty). $

\medskip

Denote by $\Omega_1$ the measurable subset of
$\Omega'_1$ such that $P_1(\Omega_1)=1$ and the convergence 2) is true for all $\omega_1\in\Omega_1$. Also by Skorokhod theorem the convergence $\Lambda_n\stackrel{d}{\to}\Lambda$ at $n\to\infty$ in
$C[0,1]$ implies the existence of probability space
$(\Omega'_2, \GA_2, \PP_2)$
 and random processes $\Gamma_n:\Omega'_2\rightarrow C[0,1]$ such that

\medskip

 $ 1) {\cal L} (\Gamma_n)={\cal L} (\Lambda_n); \ \ $

\medskip

 $ 2) \ \Gamma_n\casure \Gamma \mbox{ at } n\to\infty \mbox{ in } C[0,1]. $

\medskip

In the capacity of measurable subset  $A\subset C[0,1]$ consider the set  $\Pi[0,1]$.
Thus ${\cal L}  (\Gamma_n)(A)={\cal L} ( \Lambda_n)(A)=1$  and ${\cal L}
(\Gamma)(A)={\cal L} (\Lambda)(A)=1$, the random elements $\Gamma_n$ and
$\Gamma$ take values in $\Pi[0,1]$.
 Denote by $\Omega_2$ the measurable subset of $\Omega'_2$ such that $P_2(\Omega_2)=1$ and the convergence  2) is true for all  $\omega_2\in\Omega_2$.
 Further we will consider the probability space
 $(\Omega, \GA, \PP)$, where $\Omega=\Omega_1\times\Omega_2$, the $\sigma$-algebra $\GA$ consists of elements of $\sigma$-algebra $\GA_1\times\GA_2$
belonging to $\Omega$, the probability $\PP$ is a restriction of probability $\PP_1\otimes\PP_2$ to $\sigma$-algebra $\GA$.
 Denote

 $$G_n(t)=G'_n(\Gamma_n (t)), \ \ t\in [0,1]; $$

  $$G(t)=G'(\Gamma( t)), \ \ t\in [0,1]. $$

  Thus the trajectories of random process $\Gamma_n$ are nondecreasing, the trajectories of $G_n$ belong to the Skorokhod space $D[0,1]$.

 Let $\omega=(\omega_1,\omega_2)\in \Omega$, then the functions $g_n(t)\equiv G'_n(t)(\omega_1),$ $ g(t)\equiv G'(t)(\omega_1),$    $ \gamma_n(t)\equiv\Gamma_n(t)(\omega_2), $   $\gamma(t)\equiv \Gamma(t)(\omega_2)$ satisfy the conditions of Lemma 1.
Therefore,  $$\rho_1(g_n\circ\gamma_n,g\circ\gamma)\to 0 \mbox{ at } n\to\infty$$ and  $G_n\casure G$ at $n\to\infty$ in $D[0,1]$. This implies the convergence $X_n\stackrel{d}{\to} X$ at $n\to\infty$ in $D[0,1]$.
 
When the condition  $(B)$ is satisfied the proof is similar, but 	
instead Lemma 1 we use Lemma 2.
 The proof s completed.

{\bf NB. } The proof of Theorem 1 for the condition $(B)$ is given in \cite{silvestrov1974}, see. sections 2 and 3.

The following example shows that we can't waive the conditions (*) on the functions of space $\Pi[0,1]$.

   {\bf Example 1.} Define the sequence of the functions:
$$
 g_{n}(t)=\left\{
 \begin{array}{lr}1
\mbox{ if }& \frac{1}{2}-\frac{1}{2^n}\le t<\infty,\\
0\mbox{ if }& 0\le t<\frac{1}{2}-\frac{1}{2^n}.\\
\end{array}
\right.
$$
This sequence converges in Skorokhod metric to the function
$$
 g(t)=\left\{
 \begin{array}{lr} 1\mbox{ if }
& \frac{1}{2}\le t<\infty,\\
0\mbox{ if }& 0\le t<\frac{1}{2}.\\
\end{array}
\right.
$$
Let the sequence of functions $\gamma_n(t), \ t\in[0,1]$ is defined by equalities:  $$\gamma_{2n}(t)=\alpha_{2n}\left(\frac{1}{2}-\frac{1}{2^{2n+1}}\right)t,$$
$$\gamma_{2n+1}(t)=\alpha_{2n+1}\left(\frac{1}{2}+\frac{1}{2^{2(n+1)}}\right)t,$$ where $\alpha_n=1-\frac{1}{2^n-1}-\frac{1}{n^2}, \ n>1, \ \ \ \alpha_1=0.$ Let $\gamma(t)=\frac{1}{2}t$. It is easy to see that $\gamma_n(t)\to\gamma(t)$ at $n\to\infty$. The require $\gamma_n(1)=\gamma(1)$ is violated. In addition note that $g_{2n}(\gamma_{2n}(t))\equiv 0,$
$$
 g_{2n+1}(\gamma_{2n+1}(t))=\left\{
 \begin{array}{lr} 1 \mbox{ if }
& 2\frac{2^{2n}-1}{(2^{2n+1}+1)\alpha_{2n+1}}\le t\le 1,\\
0 \mbox{ if }& 0\le t<2\frac{2^{2n}-1}{(2^{2n+1}+1)\alpha_{2n+1}}. \\
\end{array}
\right.
$$ and the sequence $g_{2n+1}(\gamma_{2n+1}(t))\to 1$ at $n\to \infty$ in $D[0,1]$.
Since $g(\gamma(t))\equiv 1$ then $g_n\circ\gamma_n\not\to g\circ\gamma$ at $n\to\infty$ in $D[0,1]$.

The following example shows that the condition of continuity of limit function $g$ in Lemma 2 is significantly.

 {\bf Example 2. } Consider the sequence of the functions
$$
 g_{n}(t)=\left\{
 \begin{array}{lr}1\mbox{ if }
& \frac{1}{2}-\frac{1}{2^n}\le t<\infty,\\
0 \mbox{ if }& 0\le t<\frac{1}{2}-\frac{1}{2^n}.\\
\end{array}
\right.
$$
Then $g_n\stackrel{d}{\to}g$  at  $n\to\infty $ в
$D[0,\infty)$, where
$$
 g(t)=\left\{
 \begin{array}{lr} 1 \mbox{ if }
& \frac{1}{2}\le t<\infty,\\
0 \mbox{ if }& 0\le t<\frac{1}{2}.\\
\end{array}
\right.
$$
Let the sequence of the functions $\gamma_n(t), \ t\in[0,1]$ takes a constant values:
 $\gamma_n(t)=\frac{1}{2}-\frac{1}{2^{n+1}}$ and $\gamma(t)=\frac{1}{2}.$
 Then $\gamma_n(t)\to\gamma(t)$ at $n\to\infty$.
 It is easy to see that $g_n(\gamma_n(t))\equiv 0,$  $g(\gamma(t))\equiv 1$  and therefore
  $g_n\circ\gamma_n\not\to g\circ\gamma$ at $n\to\infty$  in $D[0,1]$.
 
 Consider some corollaries of this theorem.

Denote by $\pi(t)$, $t\in[0,\infty)$ the Poisson random process with intensity 1.
Consider the sequence of random processes
  $$X'_n(t)=\sum_{i=1}^{\pi(nt)}\xi_{in}, \ \ \ \ t\in[0,1], $$ where
  $\xi_{in}$ are independent identically distributed for each  $n\in {\bf N}$ random variables and $\Lambda_n$ is a sequence of random processes with nonnegative nondecreasing almost sure finite trajectories in $D[0,1]$ such that $\Lambda_n(0)=0$.
  Denote:
$$X_n(t)=\sum_{i=1}^{\pi(\Lambda_n(t))}\xi_{in}, \ \ \ \ \ t\in[0,1]. $$
   By $W'$ we will denote the Winer random process with trajectories in $D[0,\infty)$, by $W$ --- the Winer random process with random time substitution: $W(t)\equiv W'(\Lambda(t))$, $t\in[0,1]$.

  {\bf Corollary 1.}
  {\it Let  $\frac{\Lambda_n}{n}\stackrel{d}{\to}\Lambda$ at $n\to\infty$ in $D[0,1]$, $\sum_{i=1}^{n}\xi_{in}\stackrel{d}{\to}\gamma$ at $n\to\infty$, where $\gamma$  is Gaussian random variable with zero mean and  variance equal to one. Let $\Lambda_n$, $\pi$, $\xi_{in}$ and $\Lambda$, $W'$ are independent. Then
$$X_n\stackrel{d}{\to}W \mbox{  at } n\to \infty \mbox{ in } D[0,1].$$}

{\bf Proof.} The convergence
$\sum_{i=1}^{n}\xi_{in}\stackrel{d}{\to}\gamma$ at $n\to\infty$
implies (see the corollary 1 in \cite{perm}) the convergence
$X'_n\stackrel{d}{\to}W'$ at $n\to \infty$ in $D[0,\infty)$.
Since the condition $(B)$ of Theorem 1 is satisfied.

Corollary 1 can be applied to Insurance Mathematic. Consider the following model: let $n$ is a count of contracts sold by an insurance company,   $\xi_{in}$ is a extent of loss in the $i$-th insured accident for a portfolio of $n$ contracts. Assume that $\xi_{in}$ are independent identically distributed for each $n\in {\bf N}$ random values and that the count of accidents in the portfolio of $n$ contracts to the moment $t$ is described by Poisson process $\pi(nt)$ with intensity 1. Then the random process
$$X'_n(t)=\sum_{i=1}^{\pi(nt)}\xi_{in}, \ \ \ \ t\in[0,1] $$ describes the losses of insurance company, which came to moment $t$. This model does not take into account the significant for the insurance company fact that the number of insurers' contracts during the period of time
considered changes unevenly. Let $\Lambda_n(t)$ is a sequence of random processes describing the number of operating by the time  $t$ contracts. Suppose that the number of contracts increase (i. e. the trajectories of process $\Lambda_n$ are not decrease). Then the random process $$X_n(t)=\sum_{i=1}^{\pi(\Lambda_n(t))}\xi_{in}, \ \ \ \ \ t\in[0,1] $$ describes the amount of loss of the insurance company at the time $t$ by the current at that time contracts. 

Assume that $\frac{\Lambda_n}{n}\stackrel{d}{\to}\Lambda$ at $n\to\infty$, then the problem of asymptotic behavior of loss processes can be solved by means of Corollary 1. 

In terms of practical calculations the following case is interest. Suppose that with the growth of the insurance portfolio fluctuations of intensity of the process of new contracts dues become insignificant, i. e.
$\Lambda(t)= at$ where $a>0$ is a some constant. 
The Corollary 2 is a particularly case of Corollary 1.

{\bf Corollary 2.}
{\it Let  $\frac{\Lambda_n}{n}\stackrel{d}{\to}\Lambda$ at $n\to\infty$ in $D[0,1]$, $\Lambda(t)=at$, $\sum_{i=1}^{n}\xi_{in}\stackrel{d}{\to}\gamma$ at $n\to\infty$, where $\gamma$ is a Gaussian random variable  with zero mean and  variance equal to one. Let $\Lambda_n$, $\pi$ and $\xi_{in}$ are independent. Then
$$X_n\stackrel{d}{\to}\sqrt{a}W' \mbox{  at } n\to \infty$$ in $D[0,1]$ in the uniform norm.}

\end{document}